\documentclass{amsart}
\usepackage{amssymb}

\usepackage{a4}
\begin{document}
\newtheorem{theorem}{Theorem}[section]
\newtheorem{definition}[theorem]{Definition}
\newtheorem{lemma}[theorem]{Lemma}
\newtheorem{example}[theorem]{Example}
\newtheorem{remark}[theorem]{Remark}
\def\ffrac#1#2{{\textstyle\frac{#1}{#2}}}
\def\rank{\operatorname{Rank}}\def\range{\operatorname{Range}}
\def\oo{\mathfrak{o}}
\def\MM{\mathfrak{M}}
\def\AA{\mathfrak{A}}
\def\qedbox{\hbox{$\rlap{$\sqcap$}\sqcup$}}
\def\xx{x^*}\def\XX{X^*}\def\zz{z^*}\def\ZZ{Z^*}\def\tzz{\tilde z^*}\def\tZZ{\tilde Z^*}
\def\ss{s^*}\def\SS{S^*}\def\tss{\tilde s^*}\def\tSS{\tilde S^*}
\def\dss{\ddot s^*}
\makeatletter
 \renewcommand{\theequation}{%
 \thesection.\alph{equation}}
 \@addtoreset{equation}{section}
 \makeatother
\title[$k$-Curvature homogeneous manifolds]
{Complete $k$-curvature homogeneous pseudo-Riemannian manifolds $0$-modeled on an indecomposible symmetric
space}
\author{P. Gilkey and S. Nik\v cevi\'c}
\begin{address}{PG: Mathematics Department, University of Oregon,
Eugene Or 97403 USA.\newline Email: {\it gilkey@darkwing.uoregon.edu}}
\end{address}
\begin{address}{SN: Mathematical Institute, SANU,
Knez Mihailova 35, p.p. 367,
11001 Belgrade,
Serbia and Montenegro.
\newline Email: {\it stanan@mi.sanu.ac.yu}}\end{address}
\begin{abstract} For $k\ge2$, we exhibit complete $k$-curvature homogeneous neutral signature
pseudo-Riemannian manifolds which are not locally affine homogeneous (and hence not locally homogeneous). The curvature
tensor of these manifolds is modeled on that of an indecomposible symmetric space. All the local scalar Weyl curvature invariants
of these manifolds vanish. 
\end{abstract}
\maketitle
\centerline{Dedicated to Professor Sekigawa on his 60th birthday}
\section{Introduction}\label{sect-1}

\subsection{Affine manifolds} Let $\mathcal{A}:=(M,\nabla)$ be an {\it affine manifold} where $\nabla$ is a torsion free connection
on a smooth manifold
$M$. Let $\mathcal{R}_A$
be the associated curvature operator:
$$
\mathcal{R}_{\mathcal{A}}(\xi_1,\xi_2)\xi_3:=(\nabla_{\xi_1}\nabla_{\xi_2}-\nabla_{\xi_2}\nabla_{\xi_1}-\nabla_{[\xi_1,\xi_2]})\xi_3\,.
$$
Let $\nabla^i\mathcal{R}_\mathcal{A}$ be the
$i^{\operatorname{th}}$ covariant derivative of the curvature operator. If
$P\in M$, let $\nabla^i\mathcal{R}_{\mathcal{A},P}$ be the restriction of $\nabla^i\mathcal{R}_{\mathcal{A}}$ to $T_PM$. Consider
the following algebraic structure which encodes the covariant derivatives of the curvature operator up to order $k$:
$$\mathfrak{A}^k(\mathcal{A},P):=(T_PM,\mathcal{R}_{\mathcal{A},P},...,\nabla^k\mathcal{R}_{\mathcal{A},P})\,.$$
We say that $\phi:\mathfrak{A}^k(\mathcal{A}_1,P_1)\rightarrow\mathfrak{A}^k(\mathcal{A}_2,P_2)$ is an {\it affine isomorphism} if
$\phi$ is a linear map from $T_{P_1}M_1$ to
$T_{P_2}M_2$ satisfing
$$\phi^*\{\nabla_2^i\mathcal{R}_{\mathcal{A}_2,P_2}\}
=\nabla_1^i\mathcal{R}_{{\mathcal{A}_1},P_1}\quad\text{for}\quad0\le
i\le k\,.$$

\subsection{Pseudo-Riemannian manifolds} If $\mathcal{M}:=(M,g)$ is a pseudo-Riemannian manifold of signature $(p,q)$ and of
dimension $m=p+q$, let
$\nabla$ be the Levi-Civita connection, let $\mathcal{A}(\mathcal{M}):=(M,\nabla)$ be the underlying affine structure, and let
$$
R_{\mathcal{M}}(\xi_1,\xi_2,\xi_3,\xi_4):=g(\mathcal{R}_{\mathcal{A}}(\xi_1,\xi_2)\xi_3,\xi_4)
$$ 
be the curvature tensor; $R_\mathcal{M}\in\otimes^4T^*M$.
Similarly, let $\nabla^iR_{\mathcal{M}}$ be the
$i^{\operatorname{th}}$ covariant derivative of the curvature tensor. Let
$$
\mathfrak{M}^k(\mathcal{M},P):=(T_PM,g_P,R_{\mathcal{M},P},...,\nabla^kR_{\mathcal{M},P})\,.
$$
One says that
$\phi:\mathfrak{M}^k(\mathcal{M}_1,P_1)\rightarrow\mathfrak{M}^k(\mathcal{M}_2,P_2)$ is an isomorphism if $\phi$ is a linear
isomorphism from $T_{P_1}M_1$ to $T_{P_2}M_2$ so that
$$\phi^*\{g_{2,P_2}\}=g_{1,P_1}\quad\text{and}\quad
\phi^*\{\nabla_2^iR_{\mathcal{M}_2,P_2}\}=\nabla_1^iR_{\mathcal{M}_1,P_1}\quad\text{for}\quad0\le i\le k\,.$$
In this situation, the metric permits one to raise indices and conclude as well that
$$\phi^*\{\nabla_2^i\mathcal{R}_{\mathcal{M}_2,P_2}\}=\nabla_1^i\mathcal{R}_{\mathcal{M}_1,P_1}\quad\text{for}\quad0\le i\le
k\,.$$ Thus $\phi$ is also an isomorphism from $\mathfrak{A}^k(\mathcal{A}(\mathcal{M}_1),P_1)$ to
$\mathfrak{A}^k(\mathcal{A}(\mathcal{M}_2),P_2)$ of the underlying affine structure.

We shall frequently simplify the notation by setting $\mathcal{R}=\mathcal{R}_{\mathcal{A}}$ or  $R=R_{\mathcal{M}}$ when no
confusion is likely to result.

\subsection{Various notions of homogeneity} One is often interested in manifolds with a great deal of geometric symmetry. Sometimes
this symmetry arises from a transitive group action; such manifolds are called {\it homogeneous}.
\begin{definition}\label{defn-1.1}
\rm\ \begin{enumerate}\item An affine manifold $\mathcal{A}=(M,\nabla)$ is said to be  
{\it locally affine homogeneous} if given $P,Q\in M$, there is a diffeomorphism $\Phi_{P,Q}$
from a neighborhood of $P$ to a neighborhood of $Q$ so $\Phi_{P,Q}^*\nabla=\nabla$ and so $\Phi(P)=Q$.
\item A pseudo-Riemannian manifold $\mathcal{M}=(M,g)$ is said to be {\it locally homogeneous} if given $P,Q\in M$, there is a
diffeomorphism $\Phi_{P,Q}$ from a neighborhood of
$P$ to a neighborhood of $Q$ so $\Phi_{P,Q}^*g=g$ and so $\Phi(P)=Q$.
\end{enumerate}\end{definition}

There are, however, other less restrictive notions of symmetry arising from the curvature operator and
curvature tensor:

\begin{definition}\label{defn-1.2}
\ \rm\begin{enumerate}\item One says that an affine manifold $\mathcal{A}$ is  {\it
affine $k$-curvature homogeneous} if $\mathfrak{A}^k(\mathcal{A},P)$ and $\mathfrak{A}^k(\mathcal{A},Q)$ are isomorphic for
any
$P,Q\in M$.
\item One says that a pseudo-Riemannian manifold $\mathcal{M}$ is  {\it $k$-curvature homogeneous} if 
$\mathfrak{M}^k(\mathcal{M},P)$ and
$\mathfrak{M}^k(\mathcal{M},Q)$ are isomorphic for any $P,Q\in M$.
\end{enumerate}\end{definition}

One is interested finding manifolds which are affine $k$-curvature homogeneous but
not locally affine homogeneous or which are
$k$-curvature homogeneous but not locally homogeneous.

\subsection{Previous results} There are $2$-curvature
homogeneous affine manifolds which are not locally affine homogeneous \cite{GVV99,KO99,KO00,KO04,O96}. In
the Riemannian setting ($p=0$), Takagi \cite{T74} constructed 0-curvature
homogeneous complete non-compact manifolds which are not locally homogeneous; compact examples were exhibited subsequently by
Ferus, Karcher, and M\"unzer
\cite{FKM81}. Many other examples are known \cite{CMP00,KP94,KTV92,KTV92a,T97,T88,T88x,Va91}. There are no known Riemannian
manifolds which are $1$-curvature homogeneous but not locally homogeneous. In the Lorentzian setting ($p=1$)
$0$-curvature homogeneous manifolds which are not locally homogeneous were constructed by Cahen et. al. \cite{CLPT90};
$1$-curvature homogeneous manifolds which are not locally homogeneous were constructed by Bueken and Djori\'c
\cite{BD00} and by Bueken and Vanhecke \cite{BV97}.

\subsection{Curvature homogeneity and homogeneity} It is clear that local homogeneity implies $k$-curvature homogeneity for any
$k$. The following result, due to Singer
\cite{S60} in the Riemannian setting and to F. Podesta and A. Spiro
\cite{PS04} in the general context, provides a partial converse:
\begin{theorem}[Singer, Podesta-Spiro]\label{thm-1.3}
There exists an integer $k_{p,q}$ so that if $\mathcal{M}$ is a complete simply connected pseudo-Riemannian manifold of
signature $(p,q)$ which is $k_{p,q}$-curvature homogeneous, then $(M,g)$ is homogeneous.
\end{theorem}

These constants were first studied in the Riemannian setting.
Singer \cite{S60} showed $k_{0,m}<\frac12m(m-1)$; subsequently
Yamato \cite{Y89} and Gromov \cite{Gr86} established the bounds
$3m-5$ and
$\frac32m-1$ for
$k_{0,m}$, respectively. Sekigawa, Suga, and Vanhecke
\cite{SSV92,SSV95} showed any $1$-curvature homogeneous complete simply connected Riemannian manifold of dimension $m<5$ is
homogeneous; thus
$k_{0,2}=k_{0,3}=k_{0,4}=1$. We refer to the discussion in Boeckx, Vanhecke, and Kowalski
\cite{BKV96} for further details concerning $k$-curvature homogeneous manifolds in the Riemannian setting; Opozda
\cite{O97} has established an analogue of Theorem \ref{thm-1.3} in the affine setting. Observe that our definition of $k_{p,q}$
differs slightly from that given elsewhere by certain authors.

We  constructed \cite{GS04a} complete metrics of neutral signature $(p+3,p+3)$ on
$\mathbb{R}^{2p+6}$ for any $p\ge0$ which are $p+2$-curvature homogeneous but not affine $p+3$-curvature homogeneous \cite{GS04a}.
The discussion there shows $k_{p,q}\ge\min\{p,q\}$.

\subsection{Scalar invariants} One can use the metric to contract indices in pairs and form {\it scalar Weyl invariants}. Adopt
the Einstein convention and sum over repeated indices. Let $R_{i_1i_2i_3i_4}$ denote the components of the curvature
tensor. The scalar curvature
$\tau$ and the norm of the Ricci tensor
$|\rho|^2$ are given respectively by:
$$\tau=g^{i_1i_2}g^{j_1j_2}R_{i_1j_1j_2i_2}\quad\text{and}\quad|\rho|^2=g^{i_1j_1}g^{i_2j_2}g^{i_3j_3}g^{i_4j_4}R_{i_1i_2i_3j_1}
R_{i_4j_2j_3j_4}\,.$$
 There is a related result concerning scalar invariants:

\begin{theorem}[Pr\"ufer, Tricerri, and Vanhecke \cite{PTV96}]\label{thm-1.4}
If all local scalar Weyl invariants up to
order $\frac12m(m-1)$ are constant on a Riemannian manifold $\mathcal{M}$, then $\mathcal{M}$ is locally homogeneous and
$\mathcal{M}$ is determined up to local isometry by these invariants. 
\end{theorem}

This result fails in the pseudo-Riemannian setting; Koutras and McIntosh \cite{KM96} gave examples of non-flat manifolds all of
whose scalar Weyl invariants vanish; see also related examples by Pravda, Pravdov\'a, Coley, and Milson \cite{PPCM02}.

\subsection{Riemannian manifolds modeled on homogeneous spaces} One says that $\mathcal{M}$ is {\it $k$-modeled} on a
homogeneous pseudo-Riemannian manifold $\mathcal{N}$ if
$\mathfrak{M}^k(\mathcal{M},P)$ and
$\mathfrak{M}^k(\mathcal{N},Q)$ are isomorphic for any $P\in M$ and $Q\in N$; the precise $Q\in N$ being irrelevant as
$\mathcal{N}$ is homogeneous.
One has the following results in the Riemannian and Lorentzian settings:

\begin{theorem}\ \begin{enumerate}\item{\rm(Tricerri and Vanhecke \cite{TV86})}\label{thm-1.5}
If a Riemannian manifold $\mathcal{M}$ is $0$-modeled
on an irreducible symmetric space $\mathcal{S}$, then $\mathcal{M}$ is locally isometric
to $\mathcal{S}$.
\item{\rm(Cahen et al. \cite{CLPT90})} If a Lorentzian manifold $\mathcal{M}$ is $0$-modelled on an irreducible symmetric space,
then 
$\mathcal{M}$ has constant sectional curvature.
\end{enumerate}
\end{theorem}

There is a bit of technical fuss here. Recall that a pseudo-Riemannian manifold $\mathcal{M}$ is said to be {\it irreducible} if
the holonomy representation is irreducible, i.e. if $T_PM$ does not have any proper non-trivial subspace which is invariant under
the holonomy representation for any (and hence for all) $P\in M$; $\mathcal{M}$ is said to be {\it indecomposible} if there
does not exist a non-trivial decomposition of $T_PM$ which is invariant under the holonomy representation. 

These two notions are
equivalent in the Riemannian setting but are not equivalent in the higher signature setting. It is known that there are
$1$-curvature homogeneous
$3$-dimensional Lorentzian manifolds which are modeled on an indecomposible symmetric space (which is not irreducible) but which
are not locally homogeneous; see
\cite{Bu97,BD00,BV97,CLPT90} for further details.

 In this paper, we turn to the question of constructing pseudo-Riemannian manifolds which are $0$-curvature modeled on an
indecomposible symmetric space and which are $k$-curvature homogeneous for arbitrarily large $k$; our construction is motivated
by the examples described in \cite{GS04a}. We shall be defining several tensors. To simplify the discussion, we only give the
non-zero entries in these tensors up to the usual symmetries. 

\subsection{The pseudo-Riemannian manifolds $\mathcal{M}_{6+4p,f}$}  For
$p\ge1$, let
$$(x,z_0,...,z_p,\tilde z_0,...,\tilde z_p,\xx,\zz_0,...,\zz_p,\tzz_0,...,\tzz_p)$$
be coordinates on $\mathbb{R}^{6+4p}$. If $f$ is a smooth
function on
$\mathbb{R}^{p+1}$, a {\it generalized plane wave} manifold $\mathcal{M}_{6+4p,f}:=(\mathbb{R}^{6+4p},g_{6+4p,f})$ of neutral
signature $(3+2p,3+2p)$ may be defined by setting:
\begin{eqnarray*}
&&g_{6+4p,f}(\partial_x,\partial_x)=-2\{f(z_0,...,z_p)+z_0\tilde z_0+...+z_p\tilde z_p\},\quad\text{and}\\
&&g_{6+4p,f}(\partial_x,\partial_{\xx})=g_{6+4p,f}(\partial_{z_i},\partial_{\zz_i})
=g_{6+4p,f}(\partial_{\tilde z_i},\partial_{\tzz_i})=1\,.
\end{eqnarray*}

A word on notation. The {\it dual variables}
$\{\xx,\zz_i,\tzz_i\}$ enter only rather trivially; Theorem \ref{thm-1.6} below will imply that
$\nabla^iR(\cdot)$ vanishes if any entry
belongs to the span of
$\{\partial_{\xx},\partial_{\zz_i},\partial_{\tzz_i}\}$. Thus $\mathcal{M}_{6+4p,f}$  has a parallel totally isotropic
distribution of maximal dimension. The dependence of the metric on the variables
$\{\tilde z_0,...,\tilde z_p\}$ is fixed and ensures that the $0$-model space is an indecomposible symmetric space. The crucial
variables are $\{x,z_0,...,z_p\}$. 

\subsection{The geometry of the manifolds $\mathcal{M}_{6+4p,f}$}

\begin{theorem}\label{thm-1.6}\ 
\begin{enumerate} 
\item All geodesics in $\mathcal{M}_{6+4p,f}$ extend for infinite time.
\smallbreak\item $\exp_{P,\mathcal{M}_{6+4p,f}}:T_P\mathbb{R}^{6+4p}\rightarrow\mathbb{R}^{6+4p}$ is a
diffeomorphism for all $P\in\mathbb{R}^{6+4p}$.
\item The non-zero
components of
$\nabla^kR$ are:
$$\nabla^kR(\partial_x,\partial_{\xi_1},\partial_{\xi_2},\partial_x;\partial_{\xi_3},...,\partial_{\xi_{k+2}})
=-\ffrac12(\partial_{\xi_1}\cdot\cdot\cdot\partial_{\xi_{k+2}})g_{6+4p,f}(\partial_x,\partial_x)$$
for $\xi_i\in\{{z_0},...,{z_p},{\tilde z_0},...,{\tilde z_p}\}$.
\smallbreak\item All scalar Weyl invariants of $\mathcal{M}_{6+4p,f}$ vanish.
\smallbreak\item $\mathcal{M}_{6+4p,f}$ is a symmetric space if and only if $f$ is at most quadratic.
\end{enumerate}
\end{theorem}

\subsection{The symmetric space $\mathcal{S}_{6+4p}$}

\begin{theorem}\label{thm-1.7} Let $\mathcal{S}_{6+4p}:=\mathcal{M}_{6+4p,0}$ be defined by $f=0$. Then:
\ \begin{enumerate}
\item $\mathcal{S}_{6+4p}$ is an indecomposible symmetric space.
\item $\mathcal{M}_{6+4p,f}$ is $0$-modeled on
$\mathcal{S}_{6+4p}$ for any $f=f(z_0,...,z_p)$.
\end{enumerate}\end{theorem}

\subsection{The homogeneous spaces $\mathcal{H}_{6+4p,k}$} Theorems \ref{thm-1.6} and \ref{thm-1.7} show that Theorems
\ref{thm-1.4} and
\ref{thm-1.5} fail in the higher signature context. There are other interesting
properties that this family of manifolds has. Construct a sequence of pseudo-Riemannian manifolds
$\mathcal{H}_{6+4p,k}:=\mathcal{M}_{6+4p,f_k}$ by defining:
$$f_k(z_0,...,z_p):=z_1z_0^2+...+z_kz_0^{k+1}\quad\text{if}\quad1\le k\le p,$$
and as exceptional cases
\begin{eqnarray*}
&&f_{p+1}(z_0,...,z_p):=z_1z_0^2+...+z_pz_0^{p+1}+z_0^{p+3},\quad\text{and}\\
&&f_{p+2}(z_0,...,z_p):=z_1z_0^2+...+z_pz_0^{p+1}+e^{z_0}\,.
\end{eqnarray*}

The following result shows that the local isometry type of a homogeneous space need not be determined by the first few covariant
derivatives of the curvature tensor:

\begin{theorem}\label{thm-1.8}
Let $1\le k\le p+2$. Then:\begin{enumerate}
\item $\mathcal{H}_{6+4p,k}$ is $0$-modeled on the indecomposible
symmetric space
$\mathcal{S}_{6+4p}$.
\item If $j<k$, then
\begin{enumerate}\item $\mathcal{H}_{6+4p,k}$ is $j$-modeled on $\mathcal{H}_{6+4p,j}$.
\item  $\mathcal{H}_{6+4p,j}$ is not $k$-modeled on $\mathcal{H}_{6+4p,k}$.
\end{enumerate}
\item $\mathcal{H}_{6+4p,k}$ is a homogeneous space which is not symmetric.
\end{enumerate}\end{theorem}

\subsection{The manifolds $\mathcal{N}_{6+4p,\psi}$} Let $\psi\in C^\infty(\mathbb{R})$ satisfy
$$\psi^{(p+3)}(z_0)>0\quad\text{and}\quad\psi^{(p+4)}(z_0)>0\quad\text{for all}\quad z_0\in\mathbb{R}\,.$$
Let $\mathcal{N}_{6+4p,\psi}:=\mathcal{M}_{6+4p,f_\psi}$ where
$$f_\psi:=z_1z_0^2+...+z_pz_0^{p+1}+\psi(z_0)\,.$$
The following Theorem shows that
$$\alpha_{6+4p,\psi}^k(P):=\psi^{(k+p+3)}\{\psi^{(p+3)}\}^{k-1}\{\psi^{(p+4)}\}^{-k}(P)\quad\text{for}\quad k\ge2$$
forms a collection of affine invariants which determines the isometry types of these
manifolds; these invariants are not of Weyl type. Again, this does not happen in the Riemannian setting.

\begin{theorem}\label{thm-1.9}
Suppose that $\psi_i$ are real analytic for $i=1,2$ and that $\psi_i^{(p+3)}$ and $\psi_i^{(p+4)}$ are positive. The following
assertions are equivalent:
\begin{enumerate}
\item There exists a local diffeomorphism $\phi$ from $\mathcal{N}_{6+4p,\psi_1}$ to $\mathcal{N}_{6+4p,\psi_2}$
with\newline $\phi(P_1)=P_2$ and $\phi^*\nabla_{\mathcal{N}_{6+4p,\psi_2}}=\nabla_{\mathcal{N}_{6+4p,\psi_1}}$.
\item We have $\alpha_{6+4p,\psi_1}^k(P_1)=\alpha_{6+4p,\psi_2}^k(P_2)$ for $k\ge2$.
\item There exists an isometry
$\phi:\mathcal{N}_{6+4p,\psi_1}\rightarrow\mathcal{N}_{6+4p,\psi_2}$ with $\phi(P_1)=P_2$.
\end{enumerate}\end{theorem}

\subsection{Curvature and affine homogeneity} One has the following Theorem:
\begin{theorem}\label{thm-1.10}
Assume that $\psi^{(p+3)}$ and $\psi^{(p+4)}$ are positive. Then:
\begin{enumerate}
\item $\mathcal{N}_{6+4p,\psi}$ is $0$-modeled on the indecomposible symmetric space $\mathcal{S}_{6+4p}$. 
\item $\mathcal{N}_{6+4p,\psi}$ is $j$-modeled on the homogeneous space $\mathcal{H}_{6+4p,j}$ for $1\le j\le p+2$.
\item $\mathcal{N}_{6+4p,\psi}$ is $(p+2)$-curvature homogeneous.
\item The following conditions are equivalent:
\begin{enumerate}
\item $\mathcal{N}_{6+4p,\psi}$ is homogeneous.
\item $\mathcal{N}_{6+4p,\psi}$ is affine $(p+3)$-curvature homogeneous.
\item $\alpha_{6+4p,\psi}^2$ is constant.
\item $\psi^{(p+3)}=ae^{bz_0}$ for some $a,b\ne0$.
\end{enumerate}
\end{enumerate}\end{theorem}

Taking $\psi=e^{z_0}+e^{2z_0}$ constructs a manifold which is $(p+2)$-modeled on the homogeneous space
$\mathcal{N}_{6+4p,e^{z_0}}$, which is
curvature $0$-modeled on the indecomposible symmetric space $\mathcal{S}_{6+4p}$, and which is not affine $(p+3)$-curvature
homogeneous and hence not affine homogeneous.

\section{completeness}\label{sect-2}
\begin{proof}[Proof of Theorem \ref{thm-1.6}]
To simplify the notation a bit, we introduce the variables
\begin{eqnarray*}
&&s=(s_1,....,s_{2+2p}):=(z_0,...,z_p,\tilde z_0,...,\tilde z_p),\quad\text{and}\\
&&\ss=(\ss_1,...,\ss_{2+2p}):=(\zz_0,...,\zz_p,\tzz_0,...,\tzz_p)\,.
\end{eqnarray*} 
Let $1\le i\le 2+2p$. The metric then takes the form
$$
g_{6+4p,f}(\partial_x,\partial_x)=-2F(s)\quad\text{and}\quad
g_{6+4p,f}(\partial_x,\partial_{\xx})=g_{6+4p,f}(\partial_{s_i},\partial_{\ss_i})=1\,.
$$
for $F:=f(z_0,...,z_p)+z_0\tilde z_0+z_1\tilde z_1+...+z_p\tilde z_p$.  We compute the non-zero Christoffel symbols of the first
and second kinds:
\begin{eqnarray*}
&&g_{6+4p,f}(\nabla_{\partial_x}\partial_x,\partial_{s_i})=\partial_{s_i}F,\\
&&g_{6+4p,f}(\nabla_{\partial_x}\partial_{s_i},\partial_x)=
g_{6+4p,f}(\nabla_{\partial_{s_i}}\partial_x,\partial_x)=-\partial_{s_i}F,\\
&&\nabla_{\partial_x}\partial_x=\textstyle\sum_i\partial_{s_i}F\cdot\partial_{s_i^*},\quad\text{and}\\
&&\nabla_{\partial_x}\partial_{s_i}=\nabla_{\partial_{s_i}}\partial_x=- \partial_{s_i}F\cdot \partial_{\xx}\,.
\end{eqnarray*}

The curve $\gamma(t)=(x(t),s(t),\xx(t),\ss(t))$ is a geodesic if and only if
$$
0=\ddot x,\quad 0=\ddot s_i,\quad
0=\ddot x^*-2\dot x\textstyle\sum_i\dot s_i\partial_{s_i}F,\quad\text{and}\quad 0=\dss_i+\dot x\dot
x\partial_{s_i}F\,.
$$
We solve the geodesic equation with initial conditions $\gamma(0)=(\alpha,\xi,\alpha^*, \xi^*)$
and $\dot\gamma(0)=(\beta,\eta,\beta^*,\eta^*)$ by setting:
$$\begin{array}{l}
x(t)=\alpha+\beta t,\qquad s_i(t)=\xi_i+t\eta_i,\vphantom{\vrule height 11pt}\\
x^*(t)=\textstyle\alpha^*+\beta^*t+2\beta\int_{0}^t\int_{0}^\tau
    \{\textstyle\sum_i\eta_i\partial_{s_i}F(\xi+t\eta)\}d\sigma d\tau,\vphantom{\vrule height 11pt}\\
\ss_i(t)=\xi^*_i+t\eta^*_i-\textstyle\beta^2\int_{0}^t\int_{0}^\tau
    \partial_{s_i}F(\xi+t\eta)d\sigma d\tau\,.\vphantom{\vrule height 11pt}
\end{array}$$
The solution exists for all time. Furthermore, there exists a unique geodesic with
$\gamma(0)=P$ and $\gamma(1)=Q$; this establishes Assertions (1) and (2).

Since $\nabla\partial_{ x^*}=\nabla\partial_{ \ss_i}=0$, Assertion (3) follows as the quadratic terms in the Christoffel
symbols play no role in the covariant derivatives. Let
\begin{eqnarray*}
&&\mathcal{V}_1:=\operatorname{Span}\{\partial_x-\ffrac12g_{6+4p,f}(\partial_x,\partial_x)\partial_{
x^*},\partial_{s_1},...,\partial_{s_p}\},\\
&&\mathcal{V}_2:=\operatorname{Span}\{\partial_{ x^*},\partial_{ \ss_1},...,\partial_{ \ss_p}\}\,.
\end{eqnarray*}
This decomposes $\mathbb{R}^{6+4p}=\mathcal{V}_1\oplus\mathcal{V}_2$ as the direct sum of two totally isotropic subspaces.
Since $\nabla^iR$ vanishes if any entry belongs to $\mathcal{V}_2$,
$\nabla^iR$ is supported on $\mathcal{V}_1$. As $\mathcal{V}_1$ is totally isotropic, Assertion (4) follows. Assertion
(5) is immediate from Assertion (3).\end{proof}

\section{A $0$-model for $\mathcal{M}_{6+4p,f}$}\label{sect-3}

It is convenient to work in the purely algebraic setting.  Let $V$ be an $m$ dimensional vector space. Let
$$\mathfrak{M}^k:=(V,\langle\cdot,\cdot\rangle,A^0,...,A^k)$$
where $\langle\cdot,\cdot\rangle$ is a non-degenerate inner product
on $V$ and where $A^i\in\otimes^{4+i}V^*$ satisfies the appropriate symmetries of the covariant derivatives of the curvature
tensor; if $k=\infty$, then the sequence is infinite. We say that
$\mathfrak{M}$ is a $k$-model for $\mathcal{M}=(M,g)$ if for each point $P\in M$, there is an isomorphism $\phi:T_PM\rightarrow
V$ so that
$$\phi^*\langle\cdot,\cdot\rangle=g_P\quad\text{and}\quad \phi^*A^i=\nabla^iR_P\quad\text{for}
\quad 0\le i\le k\,.$$
Clearly $\mathcal{M}$ is $k$-curvature homogeneous if and only if it admits a $k$-model as one could
take $\mathfrak{M}^k:=\mathfrak{M}^k(\mathcal{M},P)$ for any $P\in M$.

\subsection{Models for the manifolds $\mathcal{M}_{6+4p,f}$}\label{sect-3.1}
Let
$$\{X,Z_0,...,Z_p,\tilde Z_0,...,\tilde Z_p,\XX,\ZZ_0,...,\ZZ_p,\tZZ_0,...,\tZZ_p\}$$ 
be a basis for $\mathbb{R}^{6+4p}$. Define a hyperbolic inner product on $\mathbb{R}^{6+4p}$ by pairing ordinary variables with
the corresponding dual variables:
$$
\langle X,\XX \rangle=
\langle{Z_i},{\ZZ_i}\rangle=\langle{\tilde Z_i},{\tZZ_i}\rangle=1\quad\text{for}\quad 0\le i\le p\,.
$$
Define an algebraic curvature tensor $A^0$ supported on $\operatorname{Span}\{X,{Z_i},{\tilde Z_i}\}$ by:
$$A^0(X,{Z_i},{\tilde Z_i},X)=1\quad\text{for}\quad0\le i\le p.
$$
Define higher order covariant derivative curvature tensors $A^i$ for $1\le i\le p$ by:
\begin{eqnarray*}
&&A^i(X,Z_0,{Z_i},X;Z_0,...,Z_0)=1,\\
&&A^i(X,Z_0,Z_0,X;{Z_i},Z_0,...,Z_0)=1,...,\\
&&A^i(X,Z_0,Z_0,X;Z_0,...,Z_0,{Z_i})=1\,.
\end{eqnarray*}
The vectors $\{Z_i,\tilde Z_i\}$ for $0\le i\le p$ are linked by $A^0$; the vectors $Z_0$ and
${Z_i}$ are linked by
$A^i$ for $1\le i\le p$. Set
\begin{eqnarray*}
&&A^{p+1}(X,Z_0,Z_0,X;Z_0,...,Z_0)=1,\quad\text{and}\\
&&A^{p+2}(X,Z_0,Z_0,X;Z_0,...,Z_0)=1\,.
\end{eqnarray*}
For $0\le k\le p+2$, we define models:
$$\mathfrak{M}^k_{6+4p}:=(\mathbb{R}^{6+4p},\langle\cdot,\cdot\rangle,A^0,...,A^k)\,.$$

\begin{proof}[Proof of Theorem \ref{thm-1.7}] Let $0\le i,j\le p$. By Theorem \ref{thm-1.6},
$$
R(\partial_x,\partial_{z_i},\partial_{\tilde z_i},\partial_x)=1\quad\text{and}\quad
R(\partial_x,\partial_{z_i},\partial_{z_j},\partial_x)=\partial_{z_i}\partial_{z_j}F
$$
where $F=f(z_0,...,z_p)+z_0\tilde z_0+...+z_p\tilde z_p$. We set
\begin{equation}\label{eqn-3.a}\begin{array}{ll}
X:=\partial_x+F\partial_{\xx},&\XX:=\partial_{\xx},\vphantom{\vrule height 12pt}\\
Z_i:=\partial_{z_i}-\frac12\textstyle\sum_j\partial_{z_i}\partial_{z_j}f\cdot\partial_{\tilde z_j},&\ZZ_i:=\partial_{\zz_i},
\vphantom{\vrule height 11pt}\\
\tilde Z_i:=\partial_{\tilde
z_i},&\tZZ_i:=\partial_{\tzz_i}+\textstyle\frac12\sum_j\partial_{z_i}\partial_{z_j}f\cdot\partial_{\zz_j}\,.
 \vphantom{\vrule height 11pt}
\end{array}\end{equation}
We show that $\mathfrak{M}^0_{6+4p}$ is a $0$-model for $\mathcal{M}_{6+4p,f}$ by noting that the non-zero components of
$g_{6+4p,f}$ and
$R$ are then given by
\begin{equation}\label{eqn-3.b}\begin{array}{l}
g_{6+4p,f}(X,\XX)=g_{6+4p,f}(Z_i,\ZZ_i)=g_{6+4p,f}(\tilde Z_i,\tZZ_i)=1,\quad\text{and}\\
R(X,Z_i,\tilde Z_i,X)=1\quad\text{for}\quad 0\le i\le p\,.\vphantom{\vrule height 11pt}
\end{array}\end{equation}

By Theorem \ref{thm-1.6}, $\mathcal{S}_{6+4p}$ is a symmetric space. As $\mathfrak{M}^0_{6+4p}$ is a $0$-model for
$\mathcal{S}_{6+4p}$ and
$\mathfrak{M}^0_{6+4p}$ is a $0$-model for $\mathcal{M}_{6+4p,f}$,
$\mathcal{S}_{6+4p}$ is a $0$-model for $\mathcal{M}_{6+4p,f}$. To show complete the proof, we must only show
$\mathfrak{M}_{6+4p}^0$ is indecomposible.

Suppose we have a non-trivial decomposition $\mathbb{R}^{6+4p}=V_1\oplus V_2$ such that
$$A^0=A_1^0\oplus
A_2^0\quad\text{and}\quad\langle\cdot,\cdot\rangle=\langle\cdot,\cdot\rangle_1\oplus\langle\cdot,\cdot\rangle_2\,.
$$
We argue for a contradiction. Denote
the natural projections induced by this decomposition by $\pi_i:\mathbb{R}^{6+4p}\rightarrow V_i$.
 Since
$$1=\langle X,X^*\rangle=\langle\pi_1X,X^*\rangle+\langle\pi_2X,X^*\rangle$$
we may assume without loss of generality
$\langle\pi_1X,X^*\rangle\ne0$. Set $\alpha:=\pi_1(X)$.
Let $\beta\in (X^*)^\perp\cap V_2$. Then $A^0(\alpha,\cdot,\beta,\alpha)=0$ as $\alpha\in V_1$ and $\beta\in V_2$. Since $\beta$
doesn't involve $X$,
\begin{eqnarray*}
&&0=A^0(\alpha,{Z_i},\beta,\alpha)=\langle\alpha,X^*\rangle^2\langle\beta,\tilde Z_i^*\rangle,\quad\text{and}\\
&&0=A^0(\alpha,\tilde Z_i,\beta,\alpha)=\langle\alpha,X^*\rangle^2\langle\beta,Z_i^*\rangle\,.
\end{eqnarray*}
Consequently $\langle\beta,\XX\rangle=0$, $\langle\beta,\ZZ_i\rangle=0$, and $\langle\beta,\tZZ_i\rangle=0$. Thus
$$\beta\in\operatorname{Span}\{\XX,\ZZ_0,...,\ZZ_p,\tZZ_0,...,\tZZ_p\}$$
so $(X^*)^\perp\cap V_2$ is totally isotropic. Since the
restriction of $\langle\cdot,\cdot\rangle$ to $V_2$ is non-degenerate and since 
$$\dim\{(X^*)^\perp\cap V_2\}\ge\dim\{V_2\}-1,$$ we conclude
that
$\dim\{V_2\}=2$. Furthermore there must exist an element of $V_2$ not in $(X^*)^\perp$. We can therefore interchange the roles of
$V_1$ and $V_2$ to see that
$\dim\{V_1\}=2$. Consequently $6+4p=\dim\{V_1\}+\dim\{V_2\}=4$ which provides the desired contradiction.
\end{proof}

Theorem \ref{thm-1.8} (1) and Theorem \ref{thm-1.10} (1) are specials cases of Theorem \ref{thm-1.7} (2).
Theorem \ref{thm-1.8} (2b) follows since $\nabla^jR_{\mathcal{H}_{6+4p,k}}=0$ if $j>k$ whereas
$\nabla^jR_{\mathcal{H}_{6+4p,j}}\ne0$. Theorem \ref{thm-1.8} (2a) and Theorem \ref{thm-1.10} (2, 3) will follow from the following
result.

\begin{lemma}\label{lem-3.1}\ \begin{enumerate}
\item If $f=\psi(z_0)+z_1z_0^2+...+z_kz_0^{k+1}$ for $1\le k\le p$, then
$\mathfrak{M}^k_{6+4p}$ is a
$k$-model for $\mathcal{M}_{6+4p,f}$.
\smallbreak\item If $f=\psi(z_0)+z_1z_0^2+...+z_pz_0^{p+1}$ and if $\psi^{(p+3)}$ is positive on $\mathbb{R}$, then
$\mathfrak{M}^{p+1}_{6+4p}$ is a
$(p+1)$-model for
$\mathcal{M}_{6+4p,f}$.
\smallbreak\item If $f=\psi(z_0)+z_1z_0^2+...+z_pz_0^{p+1}$ and if $\psi^{(p+3)}$ and $\psi^{(p+4)}$ are positive on $\mathbb{R}$,
then
$\mathfrak{M}^{p+2}_{6+4p}$ is a
$(p+2)$-model for
$\mathcal{M}_{6+4p,f}$.
\end{enumerate}\end{lemma}

\begin{proof} We adopt the notation of Equation (\ref{eqn-3.a}). The normalizations of Equation (\ref{eqn-3.b}) are then satisfied.
Suppose $f=\psi(z_0)+z_1z_0^2+...+z_kz_0^{k+1}$. If $1\le i\le k$ and $1\le j\le p$,
\begin{eqnarray*}
&&\nabla^iR(X,Z_0,Z_0,X;Z_0,...,Z_0)=\varepsilon_i,\\
&&\nabla^iR(X,Z_0,Z_j,X;Z_0,...,Z_0)=\varepsilon_{j,i},\\
&&\nabla^iR(X,Z_0,Z_0,X;Z_0,...,Z_j,...,Z_0)=\varepsilon_{j,i}
\end{eqnarray*}
where 
$\varepsilon_i=(\partial_{z_0})^{i+2}f$ and $\varepsilon_{j,i}=(\partial_{z_0})^{i+1}\partial_{z_j}f$. Note that
$$\varepsilon_{i,i}\ne0\quad\text{for}\quad 1\le i\le k\quad\text{and}\quad\varepsilon_{j,i}=0\quad\text{for}\quad1\le j<i\le
k\,.$$

To prove Assertion (1), we must define a new frame $\{{}^1X,{}^1\XX,{}^1Z_i,{}^1\tilde Z_i,{}^1\ZZ_i,{}^1\tZZ_i\}$ so that in
addition to the relations of Equation (\ref{eqn-3.b}), the only non-zero components of $\nabla^iR$ are given by
\begin{equation}\label{eqn-3.c}\begin{array}{l}
\phantom{=}\nabla^iR({}^1X,{}^1Z_0,{}^1Z_i,{}^1X;{}^1Z_0,...,{}^1Z_0)=...\\
=\nabla^iR({}^1X,{}^1Z_0,{}^1Z_0,{}^1X;{}^1Z_0,...,{}^1Z_0,{}^1Z_i)=1\,.
\vphantom{\vrule height 11pt}\end{array}\end{equation}
Set
$${}^1X:=X\quad\text{and}\quad{}^1Z_0:=Z_0+a_1Z_1+...+a_kZ_k\,.$$
To ensure $\nabla^\ell R({}^1X,{}^1Z_0,{}^1Z_0,{}^1X;{}^1Z_0,...,{}^1Z_0)=0$ for $1\le \ell\le k$, we must have:
\begin{eqnarray*}
0&=&\varepsilon_k+(k+2)\varepsilon_{k,k}a_k,\\
0&=&\varepsilon_{k-1}+(k+1)\{\varepsilon_{k,k-1}a_k+\varepsilon_{k-1,k-1}a_{k-1}\},...\\
0&=&\varepsilon_1+3\{\varepsilon_{k,1}a_k+...+\varepsilon_{1,1}a_1\}\,.
\end{eqnarray*}
Because $\varepsilon_{i,i}\ne0$ for $1\le i\le
k$, this upper triangular system of equations is recursively solvable for $a_k$, ..., $a_1$.

To ensure that
\begin{eqnarray*}
&&\nabla^iR({}^1X,{}^1Z_0,{}^1Z_i,{}^1X;{}^1Z_0,...,{}^1Z_0)=...=1,\quad\text{and}\\
&&\nabla^iR({}^1X,{}^1Z_0,{}^1Z_j,{}^1X;{}^1Z_0,...,{}^1Z_0)=...=0\quad\text{for}\quad i\ne j\,,
\end{eqnarray*}
we set ${}^1Z_i=Z_i$ for $k<i\le p$, while for $1\le i\le k$, we set
\begin{eqnarray*}
&&{}^1{Z_1}=a_{1,1}{Z_1},\quad {}^1{Z_2}=a_{2,1}{Z_1}+a_{2,2}{Z_2},\quad...\quad
{}^1{Z_k}=a_{k,1}{Z_1}+...+a_{k,k}{Z_k}\,.
\end{eqnarray*}
To ensure that ${}^1Z_k$ is properly normalized, the following relations must hold:
\begin{eqnarray*}
1&=&a_{k,k}\varepsilon_{k,k},\\
0&=&a_{k,k-1}\varepsilon_{k-1,k-1}+a_{k,k}\varepsilon_{k,k-1},...\\
0&=&a_{k,1}\varepsilon_{1,1}+...+a_{k,k}\varepsilon_{k,1}\,.
\end{eqnarray*}
This determines ${}^1Z_k$. We continue in this fashion to determine the remaining coefficients. This ensures the proper
normalizations for $\nabla^iR$ for $1\le  i\le k$.

We now return to the relations of Equation (\ref{eqn-3.b}) for $g$ and $R$. We regard $R(X,\cdot,\cdot,X)$ as defining a
neutral signature inner product on
$$\operatorname{Span}\{\partial_{z_0},...,\partial_{z_p},\partial_{\tilde z_0},...,\partial_{\tilde z_p}\}\,.$$
Since ${}^1X=\partial_x+F\partial_{\xx}$ and since $\{{}^1Z_0,...,{}^1{Z_p}\}\subset\operatorname{Span}\{Z_0,...,Z_p\}$
we may choose
$$\{{}^1\tilde Z_0,...,{}^1\tilde Z_p\}\subset\operatorname{Span}\{\tilde Z_0,...,\tilde Z_p\}$$ so the only non-zero components
of
$R$ are
$R({}^1X,{}^1Z_i,{}^1\tilde Z_i,X)=1$. Finally, we choose a dual basis 
$$\{{}^1\XX,{}^1\ZZ_0,...,{}^1\ZZ_p,{}^1\tZZ_0,...,{}^1\tZZ_p\}\subset
\text{Span}\{\XX,\ZZ_0,...,\ZZ_p,\tZZ_0,...,\tZZ_p\}$$
so the non-zero components of the metric $g$ are 
$$g({}^1X,{}^1\XX)=g({}^1Z_i,{}^1\ZZ_i)=g({}^1\tilde Z_i,{}^1\tZZ_i)=1\,.$$
Assertion (1) of the Lemma now follows.

There is a final bit of flexibility that we use in proving Assertions (2) and (3) of the Lemma. The relations of Equation
(\ref{eqn-3.b}) continue to hold. 
We rescale the basis we have constructed by setting:
$$\begin{array}{lll}
{}^2X=\varepsilon\cdot{}^1X,& {}^2\XX=\varepsilon{}^{-1}\cdot{}^1\XX,&
{}^2Z_i=\varepsilon_i\cdot{}^1Z_i,\\ 
{}^2\ZZ_i=\varepsilon_i^{-1}\cdot{}^1\ZZ_i,& 
{}^2\tilde Z_i=\varepsilon^{-2}\varepsilon_i^{-1}\cdot{}^1\tilde Z_i,& {}^2\tZZ_i=\varepsilon^2\varepsilon_i\cdot{}^1\tZZ_i\,.
\end{array}$$
The non-zero components of $g$ and of $R$ are
\begin{eqnarray*}
&&g({}^2X,{}^2\XX)=g({}^2Z_i,{}^2\ZZ_i)=g({}^2\tilde Z_i,{}^2\tZZ_i)=1,\\
&&R({}^2X,{}^2Z_i,{}^2\tilde Z_i,{}^2X)=1\,.
\end{eqnarray*}
for $0\le i\le p$.
The non-zero components of $\nabla^i R$ for
$1\le i\le p$ are
\begin{eqnarray*}
&&\nabla^iR({}^2X,{}^2Z_0,{}^2Z_i,{}^2X;{}^2Z_0,...,{}^2Z_0)=...\\
&&\qquad=\nabla^iR({}^2X,{}^2Z_0,{}^2Z_0,{}^2X;{}^2Z_0,...,{}^2Z_i)
 =\varepsilon^2\varepsilon_i\varepsilon_0^{i+1}\,.\end{eqnarray*}
The non-zero components of $\nabla^{p+1}R$ and $\nabla^{p+2}R$ are:
\begin{eqnarray*}
&&\nabla^{p+1}R(({}^2X,{}^2Z_0,{}^2Z_0,{}^2X;{}^2Z_0,...,{}^2Z_0)=\varepsilon^2\varepsilon_0^{p+3}\psi^{(p+3)},\\
&&\nabla^{p+2}R(({}^2X,{}^2Z_0,{}^2Z_0,{}^2X;{}^2Z_0,...,{}^2Z_0)=\varepsilon^2\varepsilon_0^{p+4}\psi^{(p+4)}\,.
\end{eqnarray*}
We set $\varepsilon_i:=\varepsilon^{-2}\varepsilon_0^{-i-1}$ for $1\le i\le p$ to ensure $\nabla^iR$ has the proper
normalization for $1\le i\le p$.
 Suppose that $\psi^{(p+3)}$ is positive on $\mathbb{R}$. We normalize $\nabla^{p+1}R$ and prove Assertion (2) of the Lemma by
setting: 
$$\varepsilon_0=1\quad\text{and}\quad\varepsilon=\{\psi^{(p+3)}\}^{-1/2}\,.$$
If additionally $\psi^{(p+4)}$ is positive on $\mathbb{R}$, we may set
$$\varepsilon_0:=\psi^{(p+3)}\{\psi^{(p+4)}\}^{-1}\quad\text{and}\quad
\varepsilon=\{\varepsilon_0^{p+3}\psi^{(p+3)}\}^{-1/2}$$ to
ensure that  both
$\nabla^{p+1}$ and $\nabla^{p+2}$ are normalized appropriately. This establishes Assertion (3) of the Lemma.
\end{proof}

\section{Isometries}\label{sect-4}

Let $\mathfrak{M}^\infty(\mathcal{M},P)=(T_PM,g_{\mathcal{M}},R_{\mathcal{M},P},...,\nabla^iR_{\mathcal{M},P},...)$
be the full model at a point $P$ of a pseudo-Riemannian manifold $\mathcal{M}$. This encodes complete
information about the isometry type of the manifold under certain circumstances:

\begin{lemma}\label{lem-4.1} Let $\mathcal{M}_i:=(M_i,g_i)$ be real analytic pseudo-Riemannian manifolds for $i=1,2$. Assume
there exist points $P_i\in M_i$ so
$\exp_{P_i,\mathcal{M}_i}:T_{P_i}M_i\rightarrow M_i$ is a diffeomorphism and so there exists an isomorphism $\Phi$ between
$\mathfrak{M}^\infty(\mathcal{M}_1,P_1)$ and $\mathfrak{M}^\infty(\mathcal{M}_2,P_2)$. Then
$\phi:=\exp_{P_2,\mathcal{M}_2}\circ\Phi\circ\exp_{P_1,\mathcal{M}_1}^{-1}$ is an isometry from $\mathcal{M}_1$ to
$\mathcal{M}_2$.
\end{lemma}

\begin{proof} Belger and Kowalski \cite{BeKo94} note about analytic pseudo-Riemannian metrics that the
``metric
$g$ is uniquely determined, up to local isometry, by the tensors $R$, $\nabla R$, ..., $\nabla^kR$, ... at one point.''; see also
Gray \cite{Gr73} for related work. The desired result now follows.
\end{proof}

\begin{proof}[Proof of Theorem \ref{thm-1.9}] Let $f=\psi(z_0)+z_1z_0^2+...+z_pz_0^{p+1}$. We assume $\psi^{(p+3)}$ and
$\psi^{(p+4)}$ are positive. If
$k\ge p+1$, then the non-zero components of the curvature operator $\nabla^k\mathcal{R}$ are given by
\begin{eqnarray*}
&&(\nabla_{\partial_{z_0}})^k\mathcal{R}(\partial_x,\partial_{z_0})\partial_{z_0}=
-(\nabla_{\partial_{z_0}})^k\mathcal{R}(\partial_{z_0},\partial_x)\partial_{z_0}=\psi^{(k+2)}\partial_{x^*},\quad\text{and}\\
&&(\nabla_{\partial_{z_0}})^k\mathcal{R}(\partial_x,\partial_{z_0})\partial_x=
-(\nabla_{\partial_{z_0}})^k\mathcal{R}(\partial_{z_0},\partial_x)\partial_x=-\psi^{(k+2)}\partial_{z_0^*}\,.
\end{eqnarray*}
Choose $X,Z_0\in T_P\mathbb{R}^{6+4p}$ and $\Theta\in T_P^*(\mathbb{R}^{6+4p})$ so:
\begin{equation}\label{eqn-4.b}
\Theta\{(\nabla_{Z_0})^{p+1}\mathcal{R}(X,Z_0)X\}\ne0\,.
\end{equation}
For example one could take $\Theta=dz_0^*$,
$X=\partial_x$ and $Z_0=\partial_{z_0}$. Equation  (\ref{eqn-4.b})
is an invariant of the affine $p+1$-model as it does not depend on the metric and is preserved by local affine isomorphisms.
Expand
\begin{eqnarray*}
&&X=a\partial_x+a^*\partial_{x^*}+\textstyle\sum_i\{a_i\partial_{z_i}+\tilde a_i\partial_{\tilde z_i}
   +a_i^*\partial_{\zz_i}+\tilde a_i^*\partial_{\tzz_i}\},\\
&&Z_0=b\partial_x+b^*\partial_{x^*}+\textstyle\sum_i\{b_i\partial_{z_i}+\tilde b_i\partial_{\tilde z_i}
   +b_i^*\partial_{\zz_i}+\tilde b_i^*\partial_{\tzz_i}\}\,.
\end{eqnarray*}
If $k\ge p+1$,
\begin{eqnarray*}
&&\Theta\{(\nabla_{Z_0})^k\mathcal{R}(X,Z_0)X\}=(a b_0-b a_0)b_0^k
\Theta\{(\nabla_{\partial_{z_0}})^k
  \mathcal{R}(\partial_x,\partial_{z_0})(a\partial_x+a_0\partial_{z_0}\}\\
&=&(b a_0-a b_0)b_0^k\psi^{(k+2)}\Theta(a\partial_{z_0^*}-a_0\partial_{x^*})\,.
\end{eqnarray*}
By hypothesis this is non-zero when $k=p+1$. Thus 
$$a\ne0,\quad b_0\ne0,\quad b a_0-a b_0\ne0,\quad\text{and}\quad\Theta(a\partial_{z_0^*}-a_0\partial_{x^*})\ne0\,.$$
Set $\gamma:=\Theta(a\partial_{z_0^*}-a_0\partial_{x^*})$. We may now compute:
\begin{eqnarray*}
&&\frac{\Theta\{(\nabla_{Z_0})^{k+p+1}\mathcal{R}(X,Z_0)Z_0\}\big\{\Theta\{(\nabla_{Z_0})^{p+1}\mathcal{R}(X,Z_0)Z_0\}\big\}^{k-1}}
{\big\{\Theta\{(\nabla_{Z_0})^{p+2}\mathcal{R}(X,Z_0)Z_0\}\big\}^k}\\
&=&\frac{(b a_0-a b_0)b_0^{k+p+1}\psi^{(k+p+3)}\gamma)\cdot
  \{(b a_0-a b_0)b_0^{p+1}\psi^{(p+3)}\gamma\}^{k-1}}
{\{(b a_0-a b_0)b_0^{p+2}\psi^{(p+4)}\gamma\}^k}\\
 &&\qquad=\psi^{(k+p+3)}\{\psi^{(p+3)}\}^{k-1}\{\psi^{(p+4)}\}^{-k}=\alpha_{6+4p,\psi}^k\,.
\end{eqnarray*}
This shows that $\alpha_{6+4p,\psi}^k$ is an affine invariant. Consequently Assertion (1) implies Assertion (2) in Theorem
\ref{thm-1.9}.

We now show Assertion (2) implies Assertion (3) in Theorem \ref{thm-1.9}; this will complete the proof as it is immediate that
Assertion (3) implies Assertion (1). By Lemma \ref{lem-3.1} (3), we can choose a basis $\{X,\XX,Z_i,\tilde Z_i,\ZZ_i,\tZZ_i\}$
which normalizes
$g_{6+4p,f}$ and
$\nabla^iR$ appropriately for $0\le i\le p+2$. Since 
$$\nabla^iR(X,Z_0,Z_0,X;Z_0,...,Z_0)=1\quad\text{for}\quad i=p+1,p+2,$$
we have
$$\nabla^{k+p+1}R(X,Z_0,Z_0,X;Z_0,...,Z_0)=\alpha_{6+4p,\psi}^k\quad\text{for}\quad k\ge2\,.$$
This shows that the higher covariant derivatives are controlled by $\alpha_{6+4p,\psi}^k$. Consequently if
$\alpha_{6+4p,\psi_1}^k(P_1)=\alpha_{6+4p,\psi_2}^k(P_2)$ for $k\ge2$, there is an isomorphism between
$\mathfrak{M}^\infty(\mathcal{N}_{6+4p,\psi_1},P_1)$ and $\mathfrak{M}^\infty(\mathcal{N}_{6+4p,\psi_2},P_2)$ and hence by Lemma
\ref{lem-4.1} an isometry between $(\mathcal{N}_{6+4p,\psi_1},P_1)$ and $(\mathcal{N}_{6+4p,\psi_2},P_2)$ as desired.
\end{proof}

\begin{proof}[Proof of Theorem \ref{thm-1.8} (4)] Let $f_k:=z_1z_0^2+...+z_kz_0^{k+1}$. Let $P_i\in\mathbb{R}^{6+4p}$. By Lemma
\ref{lem-4.1}, 
\begin{equation}\label{eqn-4.c}
\mathfrak{M}^i(\mathcal{M}_{6+4p,f_k},P_1)\approx\mathfrak{M}^i(\mathcal{M}_{6+4p,f_k},P_2)
\end{equation}
for $i=k$. Since $\nabla^jR=0$ for $j>k$, we may take $i=\infty$ in Equation (\ref{eqn-4.c}). Thus by Lemma
\ref{lem-4.1}, there is an isometry of $\mathcal{M}_{6+4p,f_k}$ taking $P_1$ to $P_2$. This shows
$\mathcal{M}_{6+4p,f}$ is a homogeneous space. The argument is the same if $f=z_1z_0^2+...+z_pz_0^{p+1}+z_0^{p+3}$ where we start with
$i=p+1$ in Equation (\ref{eqn-4.c}).

If $f=z_1z_0^2+...+z_pz_0^{p+1}+ae^{bz_0}$, then
$$\alpha_{6+4p,\psi}^k=b^{k+p+3}b^{(p+3)(k-1)}b^{(p+4)(-k)}$$
is independent of the point in question. We use Theorem \ref{thm-1.9} (2) to see $\mathcal{M}_{6+4p,f}$ is a homogeneous space.
\end{proof}

\begin{proof}[Proof of Theorem \ref{thm-1.10} (4)] By Theorem \ref{thm-1.9}, (4a) $\Rightarrow$ (4b)
$\Rightarrow (4c)$. Set $h=\psi^{(p+3)}$. If (4c) holds, then
$$k=\alpha_{6+4p,\psi}^2=h^{(2)}h\{h^{(1)}\}^{-2}\,.$$
We integrate the relation $h^{(2)}h=kh^{(1)}h^{(1)}$ to see there exist $(a,b)$ so
$$h(z_0)=\left\{\begin{array}{lll}ae^{bz_0}&\text{if}&k=1,\\
a(z_0+b)^{1/(1-k)}&\text{if}&k\ne1\,.\end{array}\right.$$
Since $a(z_0+b)^{1/(1-k)}$ vanishes when $z_0=-b$, these solutions are ruled out
by the assumption $h$ is always positive and smooth. Consequently
$h(z_0)=ae^{bz_0}$ and (4d) holds.  By Theorem \ref{thm-1.9}, (4d) $\Rightarrow$ (4a).
\end{proof}

\section*{Acknowledgments} Research of P. Gilkey partially supported by the
Max Planck Institute in the Mathematical Sciences (Leipzig). Research of S. Nik\v cevi\'c partially supported by MM 1646
(Srbija) and by the  DAAD  (Germany). Both authors wish to express their thanks to the Technical University of Berlin where parts
of  research reported here were conducted.

\end{document}